\documentclass{amsproc}

\usepackage{graphicx}
\usepackage[T1]{fontenc}

\newtheorem{theorem}{Theorem}

\newcommand{\U}{\mathcal U}
\newcommand{\F}{\mathcal F}
\newcommand{\Lip}{\mathrm{Lip}}
\newcommand{\phii}{\varphi}
\newcommand{\N}{\mathbb N}
\newcommand{\Z}{\mathbb Z}
\newcommand{\R}{\mathbb R}
\newcommand{\diam}{\mathrm{diam}}

\begin{document}

\title{The shark teeth is a topological IFS-attractor}

\author{Magdalena Nowak}
\address{M.Nowak:  Jan Kochanowski University, ul. \'Swi\k{e}tokrzyska 15, 25-406 Kielce, Poland \ and \ Faculty of Mathematics and Computer Science, Jagiellonian University,
ul.\L o\-jasiewicza 6, 30-348 Krak\'ow, Poland}
\email{magdalena.nowak805@gmail.com}

\author{Tomasz Szarek}
\address{Tomasz Szarek: Institute of Mathematics, University of Gda\' nsk, ul. Wita Stwosza 57, 80-952 Gda\' nsk, Poland }
\email{szarek@itl.pl}

\thanks{Tomasz Szarek has been supported by Polish Ministry of Science and
Higher Education  Grants N N201 419139. }

\subjclass[2010]{Primary 28A80; 54D05; 54F50; 54F45}
\keywords{Fractal, Iterated Function System, IFS-attractor}

\begin{abstract}
We show that the space called shark teeth is a topological IFS-attractor, that is for every open cover of $X=\bigcup_{i=1}^nf_i(X)$, its image under every suitable large composition from the family of continuous functions $\{f_1,...,f_n\}$ lies in some set from the cover. In particular, there exists a space which is not homeomorphic to any IFS-attractor but is a topological IFS-attractor.
\end{abstract}

\maketitle

Iterated function systems (IFS) are one of the most popular and simple method of constructing fractal structures, which has wide applications to data compression, computer graphics, medicine, economics, earthquake and weather prediction and many others. A compact metric space $X$ is called an {\em IFS-attractor} if $X=\bigcup_{i=1}^nf_i(X)$ for some contractions $f_1,\dots,f_n:X\to X$. In this case the family $\{f_1,\dots,f_n\}$ is~called an {\em iterated function system}. We recall that a~map $f:X\to X$ is a {\em contraction} if its Lipschitz constant
$$\Lip(f)=\sup_{x\ne y}\frac{d(f(x),f(y))}{d(x,y)}$$
is less than 1.

The notion of an iterated function system was introduced by John Hutchinson in 1981 \cite{Hutchinson} and popularized by Michael Barnsley \cite{B}. Topological properties of IFS-attractors were studied in \cite{Williams}, \cite{Hata} and \cite{BanakhNowak}. In particular the definition of topological IFS-attractor was proposed in the last paper: compact topological space $X$ is a {\em topological IFS-attractor} if $X=\bigcup_{i=1}^nf_i(X)$ for some continuous maps $f_1,\dots,f_n:X\to X$ with the property that for any open cover $\U$ of $X$ there is $m\in\N$, such that for any functions $g_1,\dots,g_m\in\{f_1,\dots,f_n\}$ the set $g_1\circ\dots\circ g_m(X)$ lies in some set $U\in\U$. 

Note that every compact, metric space $X$ is a topological IFS-attractor if for its any open cover $\U$ the diameter of the set $g_1\circ\dots\circ g_m(X)$ is less than the Lebesgue number of $\U$, for some $m\in\N$ and every $g_1,\dots,g_m\in\{f_1,\dots,f_n\}$.

It is easy to see that each IFS-attractor is a topological IFS-attractor but not  the other way around. Moreover, we show that a space called shark teeth, constructed in \cite{BanakhNowak}, which is not homeomorphic to attractor of any iterated function system is a topological IFS-attractor.

\section{The shark teeth}\label{s2}

Consider the piecewise linear periodic function
\begin{equation*}
\phii(t)=
 \begin{cases}
t-n & \text{if }t\in[n,n+\frac{1}{2}] \text{ for some }n\in\Z,  \\
n-t & \text{if }t\in[n-\frac{1}{2},n] \text{ for some }n\in\Z,
\end{cases}
\end{equation*}
whose graph looks as follows

\begin{picture}(300,80)(-30,-10)
\put(-20,0){\vector(1,0){300}}
\put(275,-10){$t$}
\put(0,-10){\vector(0,1){70}}
\put(-20,20){\line(1,-1){20}}
\put(0,0){\line(1,1){40}}
\put(40,40){\line(1,-1){40}}
\put(80,0){\line(1,1){40}}
\put(120,40){\line(1,-1){40}}
\put(160,0){\line(1,1){40}}
\put(200,40){\line(1,-1){40}}
\put(240,0){\line(1,1){20}}
\end{picture}

For every $n\in\N$ consider the function
\begin{equation*}
\varphi_n(t)=2^{-n}\varphi(2^nt),
\end{equation*}
which is a homothetic copy of the function $\varphi(t)$.

Spaces called {\em shark teeth} are constructed in \cite{BT} and are parametrized by an infinite non-decreasing sequence $(n_k)_{k=1}^\infty$. Let $I=[0,1]\times\{0\}$ be the {\em bone} of shark teeth, and for every $k\in\N$ let $M_k= \big\{\big(t,\frac{1}k\varphi_{n_k}(t)\big):t\in[0,1]\big\}$ be the $k$th {\em row} of teeth. The space shark teeth is given by the following formula
$$M =I\cup \bigcup_{k=1}^\infty M_k.$$

\begin{figure}[b]
	\includegraphics{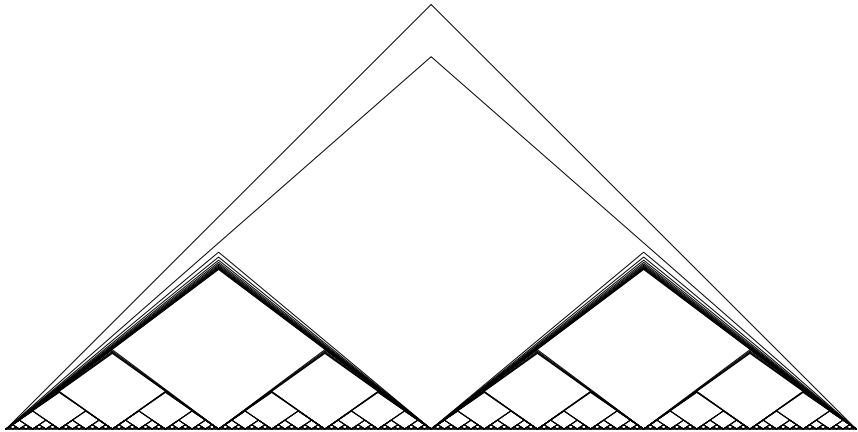}
	\caption{The space $M$}\label{shark}
\end{figure}

In \cite{BanakhNowak} is shown that the shark teeth constructed in the plane $\R^2$ with the  non-decreasing sequence
\begin{equation*}
n_k=\lfloor \log_2 \log_2 (k+1)\rfloor,\quad k\in\N,
\end{equation*}
where $\lfloor x\rfloor$ is the integer part of $x$, is  not homeomorphic to an IFS-attractor (see Figure~\ref{shark}). In other words it is not an IFS-attractor in any metric.

We show that 
\begin{theorem}
The space $M$ from \cite{BanakhNowak} is a topological IFS-attractor.
\end{theorem}

\section{Proof}

For $k\in\N$ and the sets $M_k$, $I$ and $M$, by the same names we denote the functions:
$$M_k\colon [0,1]\ni t\to \big(t,\frac{1}k\varphi_{n_k}(t)\big)\in M_k,$$
$$I\colon [0,1]\ni t\to (t,0)\in I \text{ and}$$
$$M\colon [0,1]\ni t\to I(t)\cup\bigcup_{k=1}^\infty M_k(t).$$ 

Note that for every $x\in M$ there exists unique $t_x\in[0,1]$, such that $I(t_x)=x$ or $M_k(t_x)=x$ for some $k$. Therefore we can represent every point of the space $M$ as an element from the unit interval and perhaps with positive parameter $k$. Note that for $k\neq l$ and for every $x\in M_k\cap M_l$ we have $M_k(t_x)=M_l(t_x)=I(t_x)$, because then $x$ belongs to $I$.

In three steps we will present the construction of topological IFS and prove that $M$ is its attractor.
\newline 

{\bf Step 1.} Let $\F=\{f_1,f_2,g_1,\dots,g_4,h_1,\dots,h_4\}$ be the collection of continuous functions on $M$ to itself such that for every $x\in M$ 
$$g_1|_{M\setminus M_1}(x) = M_1(0) \hspace{2cm} g_1|_{M_1}(x) = M_1\Big(\frac{\phii(t_x)}2\Big),$$
$$g_2|_{M\setminus M_1}(x) = M_1\Big(\frac{1}2\Big) \hspace{2cm} g_2|_{M_1}(x) = M_1\Big(\frac{1}2 - \frac{\phii(t_x)}2\Big),$$
$$g_3|_{M\setminus M_1}(x) = M_1\Big(\frac{1}2\Big) \hspace{2cm} g_3|_{M_1}(x) = M_1\Big(\frac{1}2 + \frac{\phii(t_x)}2\Big),$$
$$g_4|_{M\setminus M_1}(x) = M_1(1) \hspace{2cm} g_4|_{M_1}(x) = M_1\Big(1 - \frac{\phii(t_x)}2\Big).$$
Thus the union of images of $M$ under every function $g_i$ fills up the first row of the teeth $M_1 = \bigcup_{i=1}^4 g_i(M)$. Analogously we construct functions $h_i$ which fill up the second row $M_2$. Now we are going to construct functions $f_1$ and $f_2$ which cover left and right side of the rest of rows. Define $f_2(x)=f_1(x)+(\frac{1}2,0)$, so it only shifts values of function $f_1$. 

\begin{figure}
	\includegraphics[scale=0.5]{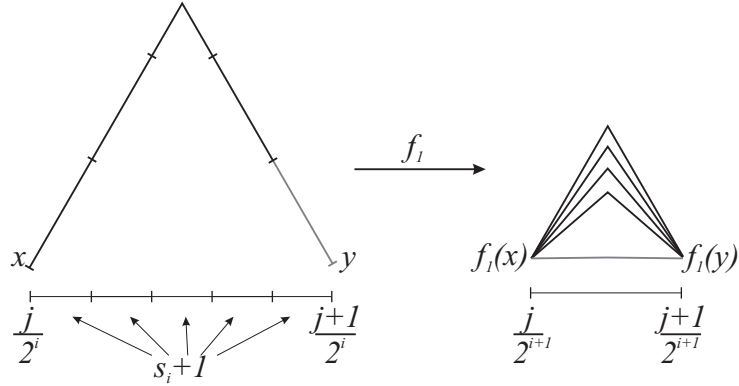}
	\caption{Tooth from $G_i$ is transform to $s_i$ teeth from $G_{i+1}$ and small part of bone $I$.}\label{tooth}
\end{figure}

For every $i\in\N$ let us define $G_i=\bigcup\{M_k: n_k=i\}$ as $i$-th {\em generation} of shark teeth. We can also treat it like a function $G_i\colon [0,1]\ni t\to \bigcup\{M_k(t): n_k=i\}\in G_i$.  Note that every row in one generation contains the same number of teeth ($2^i$). By 
$$k_i = \min\{k| n_k=i\}$$
 we denote the number of first row of teeth in $G_i$, and by 
$$N_i=|\{k| n_k=i\}|$$
 we denote the number of rows in $G_i$. Function $f_1$ has to transform every generation into the left part of next generation, so let $s_i=\frac{N_{i+1}}{N_i}$ be the number of rows from $G_{i+1}$ filled by one row from $G_i$. In our case $N_i=2^{2^{i+1}}-2^{2^i}$ and $s_i= 2^{2^{i+1}}+2^{2^i}$ for every $i\in\N$. We want the function $f_1$ to transform whole row from $G_i$ into $s_i$ rows from $G_{i+1}([0,\frac{1}2])$. Therefore, points $x,y\in M_k\cap I$ for $x\neq y$ and some positive $k$, must have distinct values $f_1(x)\neq f_1(y)$ in the same order on $I$. To obtain this, every tooth from $G_i$ must be divided into $s_i{\bf +1}$ pieces, which each of them covers one tooth from $G_{i+1}$ and the last one fills small part of bone $I$. In other words for $j=0,...,2^i-1$ a tooth from $G_i\big([\frac{j}{2^i}, \frac{j+1}{2^i}]\big)$ is transformed by $f_1$ into $s_i$ teeth from $G_{i+1}\big([\frac{j}{2^{i+1}}, \frac{j+1}{2^{i+1}}]\big)$ and bone $I\big([\frac{j}{2^{i+1}}, \frac{j+1}{2^{i+1}}]\big)$ (see Figure \ref{tooth}).

Note that for $i,j\in\N$ and for similarity $p_{i,j}(t)=\frac{t}{2^i}+\frac{j}{2^i}$, we can write $[\frac{j}{2^{i}}, \frac{j+1}{2^{i}}]=p_{i,j}([0,1])$. Moreover, define $P_{ijk}=[p_{i,j}(\frac{k}{s_i+1}),p_{i,j}(\frac{k+1}{s_i+1})]$ for $k=0,...s_i$.
Now we can present the formula for the function $f_1$:
$$f_1|_I(x)=\frac{x}2$$
and for $i\in\N$, $l=0,...,N_i-1$ and $j=0,...,2^i-1$ we have
\newline
$f_1|_{M_{k_i+l}}(x) = 
\begin{cases}
M_{k_{i+1}+ls_i+k}\Big(p_{i+1,j}\big(2\phii(\frac{s_i+1}{2}p_{i,j}^{-1}(t_x))\big)\Big), & t_x\in P_{ijk} \text{ and } k=0,...,s_i-1  \\
I(p_{i+1,j}(2\phii(\frac{s_i+1}{2}p_{i,j}^{-1}(t_x)))) & t_x\in P_{ijk} \text{ and } k=s_i
\end{cases}$

We can write that $M=\bigcup_{f\in\F}f(M)$. Indeed $\bigcup_{i=1}^4(g_i(M)\cup h_i(M)) =M_1\cup M_2$ and easy calculations can show that for every $i\in\N$ we have $f_1(G_i)=G_{i+1}([0,\frac{1}2])\cup I([0,\frac{1}2])$ and  $f_2(G_i)=G_{i+1}([\frac{1}2,1])\cup I([\frac{1}2,1])$, so 
$$f_1(M)\cup f_2(M) = \bigcup_{i=1}^\infty G_i \cup I = \overline{M \setminus (M_1\cup M_2)}.$$

{\bf Step 2.} According to the definition of functions $g_i$ and $h_i$ we have the following property for $i=0,...,4$
$$\diam~g_i(A)\leq \frac{\diam(A)}2 \text{ , } \diam~h_i(A)\leq \frac{\diam(A)}2 \text{ for every connected set} A\subset M,$$ 
so for every positive $m\in\N$ and connected set $A\subset M$ we have
\begin{equation}\label{zlozenia_g}
\diam~g_{i_1}\circ...\circ g_{i_m}(A)\leq \frac{1}{2^m}\diam(A)
\end{equation}
where $i_1,...,i_m\in\{1,...,4\}$ and analogously for functions $h_i$.

We know also the similar thing about functions $f_i$. For any positive $m\in\N$ 
\begin{equation}\label{zlozenia_f}
\diam~f_{i_1}\circ...\circ f_{i_m}(M)\leq \frac{1}{2^m}\diam(M)
\end{equation}
where $i_1,...,i_m\in\{1,2\}$. This arose due to the fact that for every $i\in\N$ and $j=0,...,2^i-1$  
$$f_1\Big(G_i\Big(\Big[\frac{j}{2^i}, \frac{j+1}{2^i}\Big]\Big)\Big) = G_{i+1}\Big(\Big[\frac{j}{2^{i+1}}, \frac{j+1}{2^{i+1}}\Big]\Big)\cup I\Big(\Big[\frac{j}{2^{i+1}}, \frac{j+1}{2^{i+1}}\Big]\Big).$$

{\bf Step 3.} Let $\U$ be an open cover of $M$. In the last step we are going to find a positive number $m$, such hat the diameter of $\phii_{i_1}\circ...\circ\phii_{i_m}(M)$ is less than the Lebesgue number $\lambda$ of $\U$, where $\phii_{i_1},...,\phii_{i_m}\in\F$. Let us consider every possible compositions of functions from $\F$. We will study the diameter of image of the space $M$ under this composition. From step 2 we know that composition of functions only from $\{g_1,...,g_4\}$, from $\{h_1,...,h_4\}$ or from $\{f_1,f_2\}$ makes half the size of the space~$M$ (see equations (\ref{zlozenia_g}) and (\ref{zlozenia_f})). Note also that for every connected set $A\subset M$ its images $g_i(A)$, $h_i(A)$ and $f_i(A)$ are contained in $\overline{M\setminus M_2}$, $\overline{M\setminus M_1}$ and $\overline{M\setminus (M_1\cup M_2)}$ respectively, so
$$\diam(g_i\circ f_j(A))=0 \hspace{1cm} \diam(g_i\circ h_j(A))=0$$
$$\diam(h_i\circ f_j(A))=0 \hspace{1cm} \diam(h_i\circ g_j(A))=0$$
because they are all singletons. This means that if the functions $g_i$, $h_i$ and $f_i$ appear in composition in the above order, the diameter of the image will be 0. It only remains for us to consider the compositions of the form $f_{i_k}\circ...\circ f_{i_1}\circ g_{j_1}\circ...\circ g_{j_n}(M)$ and analogously $f_{i_k}\circ...\circ f_{i_1}\circ h_{j_1}\circ...\circ h_{j_n}(M)$, where $i_1,...,i_k\in\{1,2\}$ and $j_1,...,j_n\in\{1,...,4\}$. Let 
$$\alpha(k)=\Lip f_1|_{G_k}= \Lip f_2|_{G_k}$$ 
be the Lipschitz constant of function $f_1$ and $f_2$ restricted to $k$-th generation. It is finite because of the definition of $f_1$. Note that the set  
$f_{i_k}\circ...\circ f_{i_1}\circ g_{j_1}\circ...\circ g_{j_n}(M)$ is contained in generation $G_{k-1}$, so we obtain
\begin{align*}
& \diam(f_{i_k}\circ...\circ f_{i_1}\circ g_{j_1}\circ...\circ g_{j_n}(M)) \leq \\
& \leq \alpha(k-1)\cdot\diam(f_{i_{k-1}}\circ...\circ f_{i_1}\circ g_{j_1}\circ...\circ g_{j_n}(M)) \leq ... \\
& \leq \alpha(k-1)\cdot...\cdot\alpha(0)\cdot\diam(g_{j_1}\circ...\circ g_{j_n}(M)) \leq \\
& \leq \prod_{i=0}^{k-1}\alpha(i)\cdot\frac{1}{2^n}\diam(M).
\end{align*}

On the other hand 
\begin{align*}
\diam(f_{i_k}\circ...\circ f_{i_1}\circ g_{j_1}\circ...\circ g_{j_n}(M))
& \leq \diam(f_{i_{k}}\circ...\circ f_{i_1}(M)) \leq \\
& \leq \frac{1}{2^k}\diam(M).
\end{align*}

Now fix $n_1\in\N$ such that $\frac{1}{2^{n_1}}\diam(M)<\lambda$ and fix $n_2\in\N$ such that $$\prod_{i=0}^{n_1-1}\alpha(i)\cdot\frac{1}{2^{n_2}}\diam(M)<\lambda.$$

Then we claim the thesis holds for $m=n_1+n_2$. Indeed, all images of $M$ under compositions only from $\{g_1,...,g_4\}$, from $\{h_1,...,h_4\}$ or from $\{f_1,f_2\}$ have diameters less than $\lambda$, because of the definition of $n_1$. Moreover $\diam(f_{i_k}\circ...\circ f_{i_1}\circ g_{j_1}\circ...\circ g_{j_{m-k}}(M))<\lambda$ for $i_1,...,i_k\in\{1,2\}$ and $j_1,...,j_n\in\{1,...,4\}$ because
\begin{enumerate}

\item if $k\leq n_1$ then 
\begin{align*}
\diam(f_{i_k}\circ...\circ f_{i_1}\circ g_{j_1}\circ...\circ g_{j_{m-k}}(M)) 
& \leq \prod_{i=0}^{k-1}\alpha(i)\cdot\frac{1}{2^{m-k}}\diam(M) \leq \\
& \leq \prod_{i=0}^{k-1}\alpha(i)\cdot\frac{1}{2^{n_2}}\diam(M)< \lambda
\end{align*}
\item  if $k>n_1$ then
\begin{align*}
\diam(f_{i_k}\circ...\circ f_{i_1}\circ g_{j_1}\circ...\circ g_{j_{m-k}}(M)) 
& \leq \frac{1}{2^{k}}\diam(M) \leq \\
& \leq \frac{1}{2^{n_1}}\diam(M)< \lambda.
\end{align*}
\end{enumerate}

Analogously we show that $\diam(f_{i_k}\circ...\circ f_{i_1}\circ h_{j_1}\circ...\circ h_{j_{m-k}}(M))<\lambda$. The others compositions transform whole space $M$ into the point so the diameter of the image of $M$ is $0<\lambda$. This ends the proof.

\section{Generalizations}\label{s3}

\begin{figure}[b]
	\includegraphics[scale=0.5]{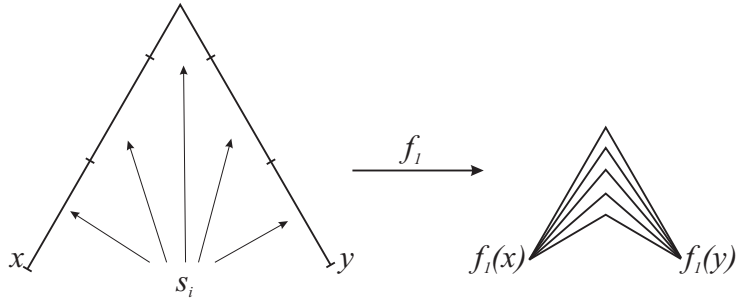}
	\caption{When $s_i$ is odd then tooth from $G_i$ is transform only to $s_i$ teeth from $G_{i+1}$.}\label{tooth2}
\end{figure}

In fact the construction above can be extended to all shark teeth.  If we try to~construct a topological IFS for shark teeth with an arbitrary sequence $(n_k)_{k=1}^\infty$, we can meet the following problems
\begin{enumerate}
\item some $G_i$ are empty.

Then we have to renumber the sequence $G_i$ such that the empty sets are~omitted. 

\item $s_i\notin\Z$.
 
Then define $s_i=\Big\lceil\frac{N_{i+1}}{N_i}\Big\rceil$, where $\lceil x \rceil$ is a minimal integer grater or equal to~$x$. Consequently, the formula for the function $f_1$ slightly changes. The~last row of teeth from every $i$-th generation has to be transformed into less than $s_i$ rows from $G_{i+1}$. It can be done by covering some rows form $G_{i+1}$ once again. 

\item $s_i$ is odd.

Then we do not have to cover a small part of bone under every tooth, so we divide every tooth from $G_i$ into $s_i$ pieces, like in the Figure \ref{tooth2}. 
\end{enumerate}

Consequently, every shark teeth is a topological IFS-attractor.

\bibliographystyle{amsplain}

\end{document}